\documentclass[11pt]{article}

\usepackage[margin=1in]{geometry}
\usepackage[justification=justified,singlelinecheck=false,font=small,labelfont=bf]{caption}
\usepackage{subcaption}
\usepackage{booktabs} 
\usepackage{tabularx}
\newcolumntype{Y}{>{\centering\arraybackslash}X} 
\newcolumntype{L}{>{\raggedright\arraybackslash}X} 
\newcolumntype{R}{>{\raggedleft\arraybackslash}X} 
\usepackage{amsmath}
\usepackage{amssymb}
\usepackage{amsfonts}
\usepackage{amsthm}
\usepackage{dsfont}
\usepackage{bold-extra}
\usepackage{cases}
\usepackage{enumerate}
\usepackage{verbatim}
\usepackage{xcolor}
\usepackage{mathtools}

\allowdisplaybreaks[4]
\usepackage{hyperref}
\hypersetup{pdfpagemode=UseNone} 
\hypersetup{pdfstartview={XYZ null null 0.85}} 

\newcommand*{\MyLaw}{\mathrm{law}}
\newcommand*{\eqlawU}{\ensuremath{\mathop{\overset{\MyLaw}{=}}}} 
\newcommand*{\eqlaw}{\mathop{\overset{\MyLaw}{\resizebox{\widthof{\eqlawU}}{\heightof{=}}{=}}}}

\newcommand*{\applaw}{\mathop{\overset{\MyLaw}{\resizebox{\widthof{\eqlawU}}{\heightof{$\approx$}}{$\approx$}}}}

\newcommand{\RR}{\mathbb{R}}

\newcommand{\bvn}{\textrm{BvN}}

\DeclareMathOperator{\E}{\mathbb{E}}
\DeclareMathOperator{\EE}{\mathbb{E}}
\DeclareMathOperator{\Prob}{\mathbb{P}}
\DeclareMathOperator{\Ind}{\mathds{1}}

\newtheorem{theorem}{Theorem}[section]
\newtheorem{proposition}[theorem]{Proposition}

\newtheorem{remark}[theorem]{Remark}

\numberwithin{equation}{section}

\begin{document}

\title{Simulation of conditional expectations under fast mean-reverting stochastic volatility models}

\author{Andrei S. Cozma\thanks{Mathematical Institute, University of Oxford, United Kingdom ({\tt andrei.s.cozma@gmail.com, christoph.reisinger@maths.ox.ac.uk})} \and Christoph Reisinger\footnotemark[1]}

\date{}
\maketitle

\begin{abstract}
In this short paper, we study the simulation of a large system of stochastic processes subject to a common driving noise and fast mean-reverting stochastic volatilities.
This model may be used to describe the firm values of a large pool of financial entities. We then seek an efficient estimator for the probability of a default, indicated by a firm value below a certain threshold, conditional on common factors.
We consider approximations where coefficients containing the fast volatility are replaced by certain ergodic averages (a type of law of large numbers), and study a correction term (of central limit theorem-type).
The accuracy of these approximations is assessed by numerical simulation of pathwise losses 
and the estimation of payoff functions as they appear in basket credit derivatives.
\end{abstract}

\section{Introduction and preliminaries}\label{sec:intro}

Consider a complete filtered probability space that is
the product of two independent probability spaces,
\[
(\Omega, \mathcal{F}, \{\mathcal{F}_t\}_{t\ge 0}, \mathbb{P}) =
(\Omega^{x,y} \times \Omega^{\dagger}, \mathcal{F}^{x,y} \otimes \mathcal{F}^{\dagger} , 
\{\mathcal{F}_t^{x,y}\otimes \mathcal{F}_t^{\dagger} \}_{t\ge 0}, \mathbb{P}^{x,y} \times \mathbb{P}^{\dagger}),
\]
such that $(\Omega^{x,y}, \mathcal{F}^{x,y}, \{\mathcal{F}^{x,y}_t\}_{t\ge 0}, \mathbb{P}^{x,y})$
supports a two-dimensional standard Brownian motion $(W^x,W^y)$ adapted to $\{\mathcal{F}_t^{x,y}\}_{t\ge 0}$ and with correlation $-1<\rho_{xy}<1$,
and $(\Omega^{\dagger}, \mathcal{F}^{\dagger}, \{\mathcal{F}^{\dagger}_t\}_{t\ge 0}, \mathbb{P}^{\dagger})$
supports an infinite i.i.d.\ sequence of two-dimensional uncorrelated standard Brownian motions $(W^{x,i}, W^{y,i})_{i\ge 1}$ adapted to 
$\{\mathcal{F}_t^{\dagger}\}_{t\ge 0}$.


For positive $N_{\hspace{-1.5pt}f} \in \mathbb{Z}$,
we study an $N_{\hspace{-1.5pt}f} \times 2$-dimensional system of SDEs of the form
\begin{equation}\label{eq_system}
\begin{aligned}
\mathrm{d}X^i_t &= \mu(V^i_t) \,\mathrm{d}t + \sigma(V^i_t)\Big( \rho_x\,\mathrm{d}W_t^x + \sqrt{1-\rho_x^2}\,\mathrm{d} W_t^{x,i} \Big),\\
\mathrm{d}V^i_t &= - \frac{k}{\epsilon} V^i_t \,\mathrm{d}t + \frac{g(V^i_t)}{\sqrt{\epsilon}} \Big( \rho_y\,\mathrm{d}W_t^y + \sqrt{1-\rho_y^2}\, \mathrm{d}W_t^{y,i}\Big),
\end{aligned}
\end{equation}
with 
$\rho_x, \rho_y \in (-1,1)$, $\epsilon, \kappa >0$ all constant; $\mu: 
\mathbb{R} \rightarrow \mathbb{R}$ and $\sigma, g: \mathbb{R} \rightarrow \mathbb{R}_+$ given functions; $((X_0^i, V_0^i))_{i\ge 1}$ are an exchangeable
infinite sequence of two-dimensional random variables that are measurable with respect to $\mathcal{F}_0 = \mathcal{F}_0^{x,y}\otimes \mathcal{F}_0^{\dagger}$.

We will consider the marginal distribution of any $X^i_t$, conditional on $\mathcal{F}^{x,y}_t$,
which is the reason for writing the Brownian driver in the decomposed way above.
Specifically, we study the setting of small $\epsilon$, a characteristic, dimensionless reversion time of $V$ to its mean.
The mean is chosen 0 here without loss of generality, but the general case is obtained by adding the constant mean to $Y$ and re-defining $\sigma$ and $\mu$.

The process $X$ is thought to describe the log-asset prices of a large portfolio of financial entities and $V$ their instantaneous stochastic volatilities. The event of $X^i$ being below a certain threshold, or barrier, $B$ models the default of that entity. Therefore, estimating marginal distributions of $X^i$ conditional on the market factors is important for the valuation and risk management of basket credit derivatives.

A simplified version of $X$ in \eqref{eq_system} with constant $\sigma$ has been considered in \cite{Bush:2011}, where an SPDE for the empirical measure in the  large pool limit is derived and used to compute tranche spreads of collateralised debt obligations, extended to jump-diffusions in \cite{Bujok:2012}.
The multilevel estimation of conditional expectations using the SDE system is analysed in \cite{Bujok:2015}, and a multilevel scheme for the SPDE in \cite{Giles:2012}.

The large pool limit under stochastic volatilities is studied in \cite{Hambly:2017}.
Computationally, this presents extra difficulties partly because of the extra dimension of the conditional expectations, but also because empirical data demonstrate a fast timescale in the volatility component (see \cite{Fouque:2000,Fouque:2003a,Fouque:2003b}), which makes accurate simulation substantially more time consuming.
Motivated by the earlier work above on ergodic limits in the context of derivative pricing (and hence parabolic PDEs), \cite{Hambly:2020} derive convergence in distribution of the conditional law of $X$ as $\epsilon\to 0$, leading to an SDE with coefficients averaged over the ergodic measure of the fast volatility process.

In this paper, we first 
present in Section \ref{sec:scheme} the simulation schemes used, including the standard Euler-Maruyama scheme
and an improved scheme which exploits exact integration of the fast process.
We then investigate in Section \ref{sec:pathwise} a number of approximations to the X process where the coefficients depending on $V$ are replaced by certain ergodic averages, and give an application to credit derivatives in Section \ref{sec:weak}.
Moreover, we compute novel correction terms, heuristically motivated by a central limit theorem-type argument, which are shown to give significantly improved results, across all scenarios considered.




For simplicity, we will restrict ourselves to the case of 
constant $g=\sqrt{2} \xi$, i.e., an Ornstein-Uhlenbeck (O--U) process $V$, and
$X_0^i =0$ and $Y_0^i=y_0$ deterministic for all $i$. 
In that case, if
we introduce a process $Z$ as the (strong) solution to 
\begin{equation} \label{OU_Z}
\mathrm{d}Z_t = -\frac{k}{\epsilon}Z_t\,\mathrm{d}t + \frac{\xi\sqrt{2}}{\sqrt{\epsilon}}\rho_y\,\mathrm{d}W_t^y,\qquad  Z_0 = 0,
\end{equation}
then $(X,Y)$ with $Y := V-Z$ satisfies, for $1\leq i\leq N_{\hspace{-1.5pt}f}$,
%
%
\begin{align}\label{eq1.1}
	\begin{dcases}
	dX^{i}_{t} \hspace{-.75em}&= \mu(Y^{i}_{t}+Z_{t}) \hspace{1pt} dt + \sigma(Y^{i}_{t}+Z_{t})\Big(\rho_{x}dW^{x}_{t}+\sqrt{1-\rho_{x}^{2}}\hspace{1pt}dW^{x,i}_{t}\Big), 
	\quad X_0^i = 0, \\
	dY^{i}_{t} \hspace{-.75em}&= -\frac{k}{\epsilon}\hspace{1pt}Y^{i}_{t}dt + \frac{\xi\sqrt{2}}{\sqrt{\epsilon}}\hspace{1pt}\sqrt{1-\rho_{y}^{2}}\hspace{1pt}dW^{y,i}_{t}, \quad Y_0^i = y_0. 
	\end{dcases}
\end{align}


Consider now the 2-dimensional empirical measure
\begin{equation}\label{eq1.2}
\nu_{N_{\hspace{-1.5pt}f},t} = \frac{1}{N_{\hspace{-1.5pt}f}}\sum_{i=1}^{N_{\hspace{-1.5pt}f}}{\delta_{X^{i}_{t},Y^{i}_{t}}}\hspace{1pt}.
\end{equation}
Using exchangeability, \cite{Hambly:2019} prove the existence of a limit 
measure
\begin{equation}\label{eq1.3}
\nu_{t} = \lim_{N_{\hspace{-1.5pt}f}\to\infty}\nu_{N_{\hspace{-1.5pt}f},t}\hspace{1pt},
\end{equation}
where the weak limit exists almost surely in $\mathbb{P}$, when $\mu=r-\sigma^2/2$ for constant $r$ and continuous bounded $\sigma$.
This follows \cite{Bush:2011} for the one-dimensional case of constant volatility, and \cite{Hambly:2017} for stochastic volatility of Cox--Ingersoll--Ross type.

Moreover, for any Borel set $A$, we have in the set-ups of \cite{Bush:2011, Ledger:2014, Hambly:2017, Hambly:2019} that
\begin{equation}\label{eq1.4}
\nu_{t}(A) = \Prob\left((X^{1}_{t},Y^{1}_{t})\in A \,|\, \mathcal{F}^{x,y}_{t}\right),
\end{equation}
where $(\mathcal{F}^{x,y}_t)_{t\ge 0}$ is here taken to be the filtration generated by the market Brownian drivers $W^{x}$ and $W^{y}$. Hence, the limit measure can be regarded as the behaviour of a single firm given the market drivers are known.

We expect these results to hold for general  $\mu$ above also, but do not provide a proof for this as it is not the focus of this paper.

%

\section{Simulation schemes for the fast O--U process}\label{sec:scheme}

Here, we first give the standard Euler--Maruyama scheme for the fast O--U processes 
and then give an alternative discretisation based on the closed-from expression for the O--U processes,
\begin{equation}
\label{eq3.1}
\begin{aligned}
Y^{i}_{t} &= y_{0}e^{-\frac{k}{\epsilon}t} + \frac{\xi\sqrt{2}}{\sqrt{\epsilon}}\hspace{1pt}\sqrt{1-\rho_{y}^{2}}\int_{0}^{t}{e^{-\frac{k}{\epsilon}(t-s)}\hspace{1pt}dW^{y,i}_{s}}, \\
Z_{t} &= \frac{\xi\sqrt{2}}{\sqrt{\epsilon}}\hspace{1pt}\rho_{y}\int_{0}^{t}{e^{-\frac{k}{\epsilon}(t-s)}\hspace{1pt}dW^{y}_{s}}.
\end{aligned}
\end{equation}


For the Euler--Maruyama scheme, we use a time mesh with timestep $\epsilon \delta t$, for some $\delta t>0$ independent of $\epsilon$.
The discrete-time approximation of $(Z_t)$ is thus generated by
\begin{equation}\label{ch4_eq_Uprocess}
\begin{aligned}
\widehat{Z}_{n} &= \widehat{Z}_{{n-1}} - k \,\delta t\,\widehat{Z}_{{n-1}} + \frac{\xi \sqrt{2}}{\sqrt{\epsilon}}\rho_y \Big(W_{t_n}^y - W_{t_{n-1}}^y\Big),
\quad n = 1,2,\ldots,
\qquad \widehat{Z}_0 = 0,
\end{aligned}
\end{equation}
where $t_n = n \delta t \epsilon$, and similar for $Y$.

The strong error is of order 1 in $\delta t$ as the diffusion coefficient is constant and the Euler--Maruyama scheme coincides with the Milstein scheme.
By choosing the time step proportionally to $\epsilon$, we found empirically that the error is asymptotically independent of $\epsilon$,
but the cost increases proportionally to $\epsilon^{-1}$.

In our second scheme, we use the closed-form expressions of $Y^{i}$ and $Z$. 
From \eqref{eq3.1} 
\begin{equation}
Y^{i}_{t} \sim \mathcal{N}\left(y_{0}e^{-\frac{k}{\epsilon}t},\hspace{1pt}\frac{\xi^{2}}{k}(1-\rho_{y}^{2})\left(1-e^{-\frac{2k}{\epsilon}t}\right)\right)
\ \stackrel{\epsilon\to 0}{\rightarrow} \ \mathcal{N}\left(0,\hspace{1pt}\frac{\xi^{2}}{k}(1-\rho_{y}^{2})\right),
\end{equation}
and
\begin{equation}\label{eq3.4}
Z_{t} \sim \mathcal{N}\left(0,\hspace{1pt}\frac{\xi^{2}}{k}\hspace{1pt}\rho_{y}^{2}\left(1-e^{-\frac{2k}{\epsilon}t}\right)\right)
\ \stackrel{\epsilon\to 0}{\rightarrow} \ \mathcal{N}\left(0,\hspace{1pt}\frac{\xi^{2}}{k}\hspace{1pt}\rho_{y}^{2}\right).
\end{equation}
Furthermore, the processes are independent across time in the limit $\epsilon\to0$ since they decorrelate exponentially fast on the time scale $\epsilon$ (see \cite{Fouque:2003a}).

For a fixed time horizon $T>0$, consider now a uniform grid $t_{n}=n\delta t$, $n\in\{0,1,\hdots,N\}$, where $T=N\delta t$. 
The discrete-time approximation processes are thus
\begin{align}\label{eq3.7}
y^{i}_{t_{n}} &= y_{0}e^{-\frac{k}{\epsilon}t_{n}} + \frac{\xi\sqrt{2}}{\sqrt{\epsilon}}\hspace{1pt}\sqrt{1-\rho_{y}^{2}}\hspace{1pt}\sum_{j=1}^{N}{e^{-\frac{k}{\epsilon}(t_{n}-t_{j-1})}\Big(W^{y,i}_{t_{j}}-W^{y,i}_{t_{j-1}}\Big)} \nonumber\\[0pt]
&= e^{-\frac{k}{\epsilon}\delta t}\bigg(y^{i}_{t_{n-1}} + \frac{\xi\sqrt{2}}{\sqrt{\epsilon}}\hspace{1pt}\sqrt{1-\rho_{y}^{2}}\Big(W^{y,i}_{t_{n}}-W^{y,i}_{t_{n-1}}\Big)\bigg),\hspace{.75em} y^{i}_{0} = y_{0}, \\
\end{align}
and
\begin{align}
\label{eq3.8}
z_{t_{n}} &= \frac{\xi\sqrt{2}}{\sqrt{\epsilon}}\hspace{1pt}\rho_{y}\sum_{j=1}^{N}{e^{-\frac{k}{\epsilon}(t_{n}-t_{j-1})}\Big(W^{y}_{t_{j}}-W^{y}_{t_{j-1}}\Big)} \nonumber\\[0pt]
&= e^{-\frac{k}{\epsilon}\delta t}\bigg(z_{t_{n-1}} + \frac{\xi\sqrt{2}}{\sqrt{\epsilon}}\hspace{1pt}\rho_{y}\Big(W^{y}_{t_{n}}-W^{y}_{t_{n-1}}\Big)\bigg),\hspace{.75em} z_{0} = 0.
\end{align}
Finally, the approximated log-asset price processes are
\begin{align}\label{eq3.9}
x^{i}_{t_{n}} &= x^{i}_{t_{n-1}} + \mu(y^{i}_{t_{n-1}}+z_{t_{n-1}})\delta t + \sigma(y^{i}_{t_{n-1}}+z_{t_{n-1}})\bigg(\rho_{x}\Big(W^{x}_{t_{n}}-W^{x}_{t_{n-1}}\Big) \nonumber\\[1pt]
&+\sqrt{1-\rho_{x}^{2}}\hspace{1pt}\Big(W^{x,i}_{t_{n}}-W^{x,i}_{t_{n-1}}\Big)\bigg),\hspace{.75em} x^{i}_{0}=0.
\end{align}

We found in experiments that if we discretize the formulae 
\eqref{eq3.1}
instead of the SDEs, this yields a lower time-discretization error.
We will therefore use the schemes \eqref{eq3.7} to \eqref{eq3.9} for the numerical tests in the subsequent sections.

\section{Pathwise conditional CDF}\label{sec:pathwise}


In this section, we give approximations to the loss function $L_T=\Prob\Big(X^{1}_{T}\leq B \,|\, \mathcal{F}^{x,y}_{T}\Big)$, i.e.\ the CDF of $X^1_T$ conditional on the market factors $W^x$ and $W^y$, using ergodic averages of coefficients and a correction term from a central limit theorem.

\subsection{Conditional CDF and Monte Carlo estimators}\label{subsec:CDF}

Let $B\in\RR$ and consider 
the loss function at time $T$ for a default level $B$,
\begin{equation}\label{eq4.2}
L_{N_{\hspace{-1.5pt}f},T} = \frac{1}{N_{\hspace{-1.5pt}f}}\sum_{i=1}^{N_{\hspace{-1.5pt}f}}{\Ind_{X^{i}_{T}\leq B}}\hspace{1pt},
\end{equation}
i.e., the proportion of companies that are in default at time $T$. Since $\big\{\hspace{-1pt}\Ind_{X^{i}_{T}\leq B} : 1\leq i\leq N_{\hspace{-1.5pt}f}\big\}$ are conditionally (on $\mathcal{F}^{x,y}_{T}$) independent and identically distributed random variables, Birkhoff's Ergodic Theorem 
(see \cite[Section V.3]{Shiryaev:96})
 implies that  the limiting loss function can be regarded as the marginal CDF, i.e.,
\begin{equation}\label{eq4.3}
L_{T} = \lim_{N_{\hspace{-1.5pt}f}\to\infty} L_{N_{\hspace{-1.5pt}f},T} = \Prob\Big(X^{1}_{T}\leq B \,|\, \mathcal{F}^{x,y}_{T}\Big).
\end{equation}
We use a conditional Monte Carlo technique to estimate the marginal CDF. Denote by $\mathcal{F}^{x,y,y_{1}}$ the filtration generated by the Brownian motions $W^{x}$, $W^{y}$, $W^{y,1}$. Then
\begin{align}\label{eq4.4}
\Prob\Big(X^{1}_{T}\leq B \,|\, \mathcal{F}^{x,y}_{T}\Big) &= \EE\left[\EE\left[\Ind_{X^{1}_{T}\leq B} \,|\, \mathcal{F}^{x,y,y_{1}}_{T}\right]\big|\, \mathcal{F}^{x,y}_{T}\right] \nonumber\\[2pt]
&= \EE\left[\Prob\Big(X^{1}_{T}\leq B \,|\, \mathcal{F}^{x,y,y_{1}}_{T}\Big)\big|\, \mathcal{F}^{x,y}_{T}\right].
\end{align}
Conditional on the $\sigma$-algebra $\mathcal{F}^{x,y,y_{1}}_{T}$, noting $W^{x,1}$ independent of $W^{y,1}$ and $W^y$, 
\begin{equation}\label{eq4.5}
\int_{0}^{T}{\sigma(Y^{1}_{t}+Z_{t})\hspace{1pt}dW^{x,1}_{t}} \,\eqlaw\, \sqrt{\int_{0}^{T}{\sigma^{2}(Y^{1}_{t}+Z_{t})\hspace{1pt}dt}}\,W_{1},
\end{equation}
where $W_{1}$ is a standard normal random variable. Hence, we deduce from \eqref{eq1.1} that
\begin{equation}\label{eq4.6}
\Prob\Big(X^{1}_{T}\leq B \,|\, \mathcal{F}^{x,y,y_{1}}_{T}\Big) = \Phi\left(\frac{B-\int_{0}^{T}{\mu(Y^{1}_{t}+Z_{t})\hspace{1pt}dt}-\rho_{x}\int_{0}^{T}{\sigma(Y^{1}_{t}+Z_{t})\hspace{1pt}dW^{x}_{t}}}{\sqrt{(1-\rho_{x}^{2})\int_{0}^{T}{\sigma^{2}(Y^{1}_{t}+Z_{t})\hspace{1pt}dt}}}\right),
\end{equation}
where $\Phi$ is the standard normal CDF. Using the discretizations from \eqref{eq3.7} and \eqref{eq3.8},
\begin{align}\label{eq4.7}
\Prob\Big(X^{1}_{T}\leq B \,|\, \mathcal{F}^{x,y,y_{1}}_{T}\Big) &\approx& \\
&&\hspace{-2.3cm} \Phi\left(\frac{B-\delta t\sum_{n=0}^{N-1}{\mu(y^{1}_{t_{n}}+z_{t_{n}})}-\rho_{x}\sum_{n=0}^{N-1}{\sigma(y^{1}_{t_{n}}+z_{t_{n}})\hspace{1pt}\Big(W^{x}_{t_{n+1}}-W^{x}_{t_{n}}\Big)}}{\sqrt{(1-\rho_{x}^{2})\delta t\sum_{n=0}^{N-1}{\sigma^{2}(y^{1}_{t_{n}}+z_{t_{n}})}}}\right).
\end{align}
The marginal CDF, i.e., the outer expectation in \eqref{eq4.4}, is estimated by a Monte Carlo average over a sufficiently large number of samples of $W^{y,1}$. As an aside, we can estimate the marginal density function by differentiating \eqref{eq4.7} with respect to $B$.


\subsection{Ergodic averages}\label{subsec:ergCDF}

We will define approximations to the process by averaging SDE coefficients over the ergodic distribution of the O--U process,
\begin{equation}\label{eq4.7.1}
\langle f \rangle_{Y} = \int_{-\infty}^{\infty}{f(y)\phi_{Y\hspace{-1pt}}(y)\hspace{1pt}dy},
\end{equation}
where
$\phi_{Y\hspace{-1pt}}$ is the centered normal density with variance $\xi^{2}(1-\rho_{y}^{2})/k$.

\textit{Linear $Y$-average.} We first approximate the marginal CDF (in $x$) by using an ergodic $Y^{1}$ average (abbreviated \emph{erg$_1$\!Y}) over its stationary distribution, namely
\begin{equation}\label{eq4.8}
\int_{0}^{T}{\sigma(Y^{1}_{t}+Z_{t})\hspace{1pt}dW^{x}_{t}} \approx \int_{0}^{T}{\langle\sigma(\cdot+Z_{t})\rangle_{Y}\hspace{1pt}dW^{x}_{t}},
\end{equation}
which matches the first conditional (on $\mathcal{F}^{x,y}_{T}$) moment of the stochastic integral in the limit $\epsilon\to0$. Hence, we obtain
\begin{align}\label{eq4.10}
\Prob\Big(X^{1}_{T}\leq B \,|\, \mathcal{F}^{x,y}_{T}\Big) &\approx \Phi\left(\frac{B-\int_{0}^{T}{\langle\mu(\cdot+Z_{t})\rangle_{Y}\hspace{1pt}dt}-\rho_{x}\int_{0}^{T}{\langle\sigma(\cdot+Z_{t})\rangle_{Y}\hspace{1pt}dW^{x}_{t}}}{\sqrt{(1-\rho_{x}^{2})\int_{0}^{T}{\langle\sigma^{2}(\cdot+Z_{t})\rangle_{Y}\hspace{1pt}dt}}}\right) \nonumber\\[2pt]
&\hspace{-2.8 cm} \approx \Phi\left(\frac{B-\delta t\sum_{n=0}^{N-1}{\langle\mu(\cdot+z_{t_{n}})\rangle_{Y}}-\rho_{x}\sum_{n=0}^{N-1}{\langle\sigma(\cdot+z_{t_{n}})\rangle_{Y}\hspace{1pt}\Big(W^{x}_{t_{n+1}}-W^{x}_{t_{n}}\Big)}}{\sqrt{(1-\rho_{x}^{2})\delta t\sum_{n=0}^{N-1}{\langle\sigma^{2}(\cdot+z_{t_{n}})\rangle_{Y}}}}\right).
\end{align}
\textit{Quadratic $Y$-average.}  Alternatively, we will use a quadratic ergodic $Y^{1}$ average (abbreviated \emph{erg$_2$\!Y}), namely
\begin{equation}\label{eq4.10.1}
\int_{0}^{T}{\sigma(Y^{1}_{t}+Z_{t})\hspace{1pt}dW^{x}_{t}} \approx \int_{0}^{T}{\langle\sigma^{2}(\cdot+Z_{t})\rangle_{Y}^{\frac{1}{2}}\hspace{1pt}dW^{x}_{t}},
\end{equation}
which matches the first and second unconditional moments of the stochastic integral in the limit $\epsilon\to0$, to obtain
\begin{align}\label{eq4.11}
\Prob\Big(X^{1}_{T}\leq B \,|\, \mathcal{F}^{x,y}_{T}\Big) &\approx \Phi\left(\frac{B-\int_{0}^{T}{\langle\mu(\cdot+Z_{t})\rangle_{Y}\hspace{1pt}dt}-\rho_{x}\int_{0}^{T}{\langle\sigma^{2}(\cdot+Z_{t})\rangle^{\frac{1}{2}}_{Y}\hspace{1pt}dW^{x}_{t}}}{\sqrt{(1-\rho_{x}^{2})\int_{0}^{T}{\langle\sigma^{2}(\cdot+Z_{t})\rangle_{Y}\hspace{1pt}dt}}}\right) \nonumber\\[2pt]
&\hspace{-2.8 cm} \approx \Phi\left(\frac{B-\delta t\sum_{n=0}^{N-1}{\langle\mu(\cdot+z_{t_{n}})\rangle_{Y}}-\rho_{x}\sum_{n=0}^{N-1}{\langle\sigma^{2}(\cdot+z_{t_{n}})\rangle^{\frac{1}{2}}_{Y}\hspace{1pt}\Big(W^{x}_{t_{n+1}}-W^{x}_{t_{n}}\Big)}}{\sqrt{(1-\rho_{x}^{2})\delta t\sum_{n=0}^{N-1}{\langle\sigma^{2}(\cdot+z_{t_{n}})\rangle_{Y}}}}\right).
\end{align}

\textit{Linear $Y$ and $Z$-average.} Third, we approximate the marginal CDF by using an ergodic $Y^{1}$ and $Z$ average (abbreviated \emph{erg$_1$\!YZ}) over their stationary distribution,
\begin{equation}\label{eq4.12}
\bar{f} = \langle\langle f(\cdot+Z) \rangle_{Y}\rangle_{Z} = \langle f \rangle_{Y\hspace{-1pt}+Z} = \int_{-\infty}^{\infty}{f(y)\phi_{Y\hspace{-1pt}+Z}(y)\hspace{1pt}dy},
\end{equation}
where
$\phi_{Y\hspace{-1pt}+Z}$ is the centered normal density with variance $\xi^{2}/k$.
Hence, we obtain
\begin{equation}\label{eq4.14}
\Prob\Big(X^{1}_{T}\leq B \,|\, \mathcal{F}^{x,y}_{T}\Big) \approx \Phi\left(\frac{B-\bar{\mu}T-\rho_{x}\bar{\sigma}W^{x}_{T}}{\sqrt{(1-\rho_{x}^{2})\overline{\sigma^{2}}T}}\right).
\end{equation}
\textit{Quadratic $Y$ and $Z$-average.}  Alternatively, we will use a quadratic ergodic $Y^{1}$ and $Z$ average (abbreviated \emph{erg$_2$\!YZ}) in the stochastic integral to obtain
\begin{equation}\label{eq4.15}
\Prob\Big(X^{1}_{T}\leq B \,|\, \mathcal{F}^{x,y}_{T}\Big) \approx \Phi\left(\frac{B-\bar{\mu}T-\rho_{x}\overline{\sigma^{2}}^{\frac{1}{2}}W^{x}_{T}}{\sqrt{(1-\rho_{x}^{2})\overline{\sigma^{2}}T}}\right).
\end{equation}


\subsection{Approximation of marginal CDF by a CLT-type argument}\label{subsec:appCDF}

Here, we introduce an approximation to the marginal CDF in $Y^{1}$ (abbreviated \emph{appY}), and hence to the limiting loss function,
by adding a correction term from a central limit theorem (CLT).
We note that,
as $\epsilon\to0$, the process $(Y^{1}_{t})_{0\leq t\leq T}$ decorrelates exponentially fast, on the time scale $\epsilon$.
Arguing informally with the central limit theorem under strong mixing (see, e.g., \cite[Theorem 27.5]{billingsley}),
we approximate for small $\epsilon$, conditional on $\mathcal{F}^{x,y}_{T}$,
\begin{equation}\label{eq4.7.6}
\frac{\int_{0}^{T}{\sigma(Y^{1}_{t}+Z_{t})\hspace{1pt}dW^{x}_{t}} - \int_{0}^{T}{\langle\sigma(\cdot+Z_{t})\rangle_{Y}\hspace{1pt}dW^{x}_{t}}
}{
\sqrt{\int_{0}^{T}{\Big(\langle\sigma^{2}(\cdot+Z_{t})\rangle_{Y}-\langle\sigma(\cdot+Z_{t})\rangle_{Y}^{2}\Big)dt}}}
\,\applaw\,
W_{1},
\end{equation}
where $W_{1}$ is a standard normal random variable. Similarly, for small $\epsilon$ and conditional on $\mathcal{F}^{x,y}_{T}$, we use
\begin{equation}\label{eq4.7.7}
\int_{0}^{T}{f(Y^{1}_{t}+Z_{t})\hspace{1pt}dt} \,\applaw\, \int_{0}^{T}{\langle f(\cdot+Z_{t})\rangle_{Y}\hspace{1pt}dt}.
\end{equation}
Note that, for any $c_{0},c_{1}\in\RR$,
\begin{equation}\label{eq4.7.8}
\EE\Big[\Phi\big(c_{0}-c_{1}W_{1}\big)\Big] = \Phi\left(\frac{c_{0}}{\sqrt{1+c_{1}^{2}}}\right).
\end{equation}
Combining \eqref{eq4.4}, \eqref{eq4.6} and \eqref{eq4.7.6}--\eqref{eq4.7.8}
yields
\begin{align}\label{eq4.7.9}
\Prob\Big(X^{1}_{T}\leq B \,|\, \mathcal{F}^{x,y}_{T}\Big) &\approx \Phi\left(\frac{B-\int_{0}^{T}{\langle\mu(\cdot+Z_{t})\rangle_{Y}\hspace{1pt}dt}-\rho_{x}\int_{0}^{T}{\langle\sigma(\cdot+Z_{t})\rangle_{Y}\hspace{1pt}dW^{x}_{t}}}{\sqrt{\int_{0}^{T}{\langle\sigma^{2}(\cdot+Z_{t})\rangle_{Y}\hspace{1pt}dt}-\rho_{x}^{2}\int_{0}^{T}{\langle\sigma(\cdot+Z_{t})\rangle_{Y}^{2}\hspace{1pt}dt}}}\right) \nonumber\\[2pt]
&\hspace{-3 cm} \approx \Phi\left(\frac{B-\delta t\sum_{n=0}^{N-1}{\langle\mu(\cdot+z_{t_{n}})\rangle_{Y}}-\rho_{x}\sum_{n=0}^{N-1}{\langle\sigma(\cdot+z_{t_{n}})\rangle_{Y}\hspace{1pt}\Big(W^{x}_{t_{n+1}}-W^{x}_{t_{n}}\Big)}}{\sqrt{\delta t\sum_{n=0}^{N-1}{\langle\sigma^{2}(\cdot+z_{t_{n}})\rangle_{Y}}-\rho_{x}^{2}\delta t\sum_{n=0}^{N-1}{\langle\sigma(\cdot+z_{t_{n}})\rangle_{Y}^{2}}}}\right).
\end{align}

\subsection{Exponential Ornstein--Uhlenbeck model}\label{subsec:expOU}

Henceforth, we consider an exponential Ornstein--Uhlenbeck stochastic volatility model for the dynamics of the asset price processes. The drift coefficient is
$\mu(y) = -\frac{1}{2}\hspace{1pt}\sigma^{2}(y)$,
whereas the diffusion coefficient is
$\sigma(y) = me^{y}$,
see \cite{Masoliver:2006}. We do not have a closed-form formula for the conditional CDF is not available under this model. 

We substitute the specific coefficients into the above formulae for the conditional CDF and use moment generating functions. From \eqref{eq4.7}, we find an estimate for the conditional (on $\mathcal{F}^{y_{1}}_{T}$) marginal (in $x$) CDF,
\begin{align}\label{eq4.18}
\Prob\Big(X^{1}_{T}\leq B \,|\, \mathcal{F}^{x,y,y_{1}}_{T}\Big) &\approx \Phi\left(\frac{Bm^{-1}+\frac{1}{2}\hspace{1pt}m\delta t\sum_{n=0}^{N-1}{e^{2y^{1}_{t_{n}}+2z_{t_{n}}}}}{\sqrt{(1-\rho_{x}^{2})\delta t\sum_{n=0}^{N-1}{e^{2y^{1}_{t_{n}}+2z_{t_{n}}}}}}\right. \nonumber\\[2pt]
&\left.-\frac{\rho_{x}}{\sqrt{1-\rho_{x}^{2}}}\hspace{1pt}\frac{\sum_{n=0}^{N-1}{e^{y^{1}_{t_{n}}+z_{t_{n}}}\hspace{1pt}\Big(W^{x}_{t_{n+1}}-W^{x}_{t_{n}}\Big)}}{\sqrt{\delta t\sum_{n=0}^{N-1}{e^{2y^{1}_{t_{n}}+2z_{t_{n}}}}}}\right).
\end{align}
From \eqref{eq4.7.9}, we find an estimate for the approximate conditional CDF,
\begin{align}\label{eq4.18.1}
\Prob\Big(X^{1}_{T}\leq B \,|\, \mathcal{F}^{x,y}_{T}\Big) &\approx \Phi\left(\frac{Bm^{-1}e^{-\frac{\xi^{2}}{k}(1-\rho_{y}^{2})}+\frac{1}{2}\hspace{1pt}me^{\frac{\xi^{2}}{k}(1-\rho_{y}^{2})}\delta t\sum_{n=0}^{N-1}{e^{2z_{t_{n}}}}}{\sqrt{\left(1-\rho_{x}^{2}e^{-\frac{\xi^{2}}{k}(1-\rho_{y}^{2})}\right)\delta t\sum_{n=0}^{N-1}{e^{2z_{t_{n}}}}}}\right. \nonumber\\[2pt]
&\left.-\frac{\rho_{x}}{\sqrt{e^{\frac{\xi^{2}}{k}(1-\rho_{y}^{2})}-\rho_{x}^{2}}}\hspace{1pt}\frac{\sum_{n=0}^{N-1}{e^{z_{t_{n}}}\hspace{1pt}\Big(W^{x}_{t_{n+1}}-W^{x}_{t_{n}}\Big)}}{\sqrt{\delta t\sum_{n=0}^{N-1}{e^{2z_{t_{n}}}}}}\right).
\end{align}
\textit{$Y$-averages.} From \eqref{eq4.10} and \eqref{eq4.11}, we find an estimate for the conditional CDF with the (linear and quadratic) ergodic $Y^{1}$ average,
\begin{align}\label{eq4.19}
\Prob\Big(X^{1}_{T}\leq B \,|\, \mathcal{F}^{x,y}_{T}\Big) &\approx \Phi\left(\frac{Bm^{-1}e^{-\frac{\xi^{2}}{k}(1-\rho_{y}^{2})}+\frac{1}{2}\hspace{1pt}me^{\frac{\xi^{2}}{k}(1-\rho_{y}^{2})}\delta t\sum_{n=0}^{N-1}{e^{2z_{t_{n}}}}}{\sqrt{(1-\rho_{x}^{2})\delta t\sum_{n=0}^{N-1}{e^{2z_{t_{n}}}}}}\right. \nonumber\\[2pt]
&\left.-\frac{\rho_{x}}{\sqrt{1-\rho_{x}^{2}}}\hspace{1pt}\frac{e^{-\lambda\frac{\xi^{2}}{2k}(1-\rho_{y}^{2})}\sum_{n=0}^{N-1}{e^{z_{t_{n}}}\hspace{1pt}\Big(W^{x}_{t_{n+1}}-W^{x}_{t_{n}}\Big)}}{\sqrt{\delta t\sum_{n=0}^{N-1}{e^{2z_{t_{n}}}}}}\right),
\end{align}
where $\lambda=0$ for the quadratic average and $\lambda=1$ for the linear average.

\textit{$Y$ and $Z$-averages.} Finally, from \eqref{eq4.14} and \eqref{eq4.15}, we find an estimate for the marginal CDF with the (linear or quadratic) ergodic $Y^{1}$ and $Z$ average,
\begin{equation}\label{eq4.20}
\Prob\Big(X^{1}_{T}\leq B \,|\, \mathcal{F}^{x,y}_{T}\Big) \approx \Phi\left(\frac{Bm^{-1}e^{-\frac{\xi^{2}}{k}}+\frac{1}{2}\hspace{1pt}me^{\frac{\xi^{2}}{k}}T}{\sqrt{(1-\rho_{x}^{2})T}} - \frac{\rho_{x}}{\sqrt{1-\rho_{x}^{2}}}\hspace{1pt}\frac{e^{-\lambda\frac{\xi^{2}}{2k}}W^{x}_{T}}{\sqrt{T}}\right).
\end{equation}
The formulae \eqref{eq4.18}--\eqref{eq4.20} indicate that the approximation errors as well as the difference between the approximations will increase with $|\rho_{x}|$.

\subsection{Pathwise numerical tests}\label{subsec:numericsStrong}

A motivation for considering pathwise tests of the different approximations is the filtering interpretation of the equations.

We fix the time horizon $T=1$ and the default level $B=-0.1$, and assign the following values to the underlying model parameters:
\begin{equation}\label{eq4.21}
y_{0}=0.2,\; m=0.1,\; k=1.0,\; \xi=0.26,\; \rho_{x}=0.9,\; \rho_{y}=0.5,\; \rho_{xy}=-0.6;
\end{equation}
we vary $\epsilon$. We refer to \cite{Fouque:2000,Fouque:2003b} for data that suggest a mean-reversion time of a few days for the S\&P500.

To produce the results in Table \ref{table1}, we fixed the paths for $(W^x,W^y)$, generated by standard sampling of i.i.d.\ normal increments,
and then produced $4 \cdot 10^5$ samples of $W^{y,1}$ to estimate the outer expectation in \eqref{eq4.4}, using the time stepping approximation \eqref{eq4.7}.
The number of samples was chosen such that the relative statistical error, estimated as the corrected sample standard deviation of the estimator divided by the value
itself, was below $0.15\%$.

Our tests with different $\epsilon$ suggest that the number of time steps should scale with $\epsilon^{-1}$ for uniform accuracy.
More specifically, for a fraction $\epsilon$ of a year, $40$ time steps were required for a sufficiently small time-discretization error that matches the statistical error.

The computations were carried out in MATLAB R2016b on a laptop with the following specifications: Intel(R) Core(TM) i7-6700HQ CPU 2.60GHz, 8GB RAM, running Windows 10 (64 bit).
The computations below took several hours to compute the `true' loss, which is why we considered only three `outer' sample paths.
The computation time for the various approximations was negligible as no inner sampling was required. This gain in efficiency is a major motivation for
the approximations studied in this paper.

\begin{figure}[htb]
\vspace{-3.5 cm}
\begin{center}
\includegraphics[scale=0.5]{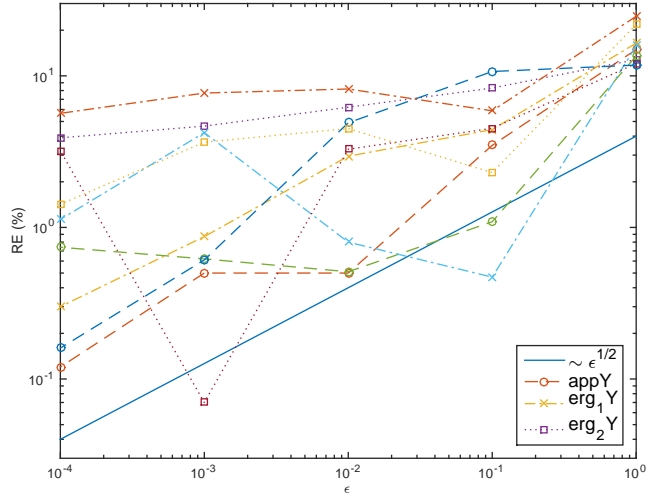}
\end{center}
\vspace{-3.5 cm}
%
%
\caption{Double logarithmic plot of the relative errors (RE) of three approximations to the marginal CDF, 
for three realisations of $(W^x,W^y)$.
Shown are the errors of the CLT approximation, appY (circles), and the linear and quadratic $Y$-averages, erg$_1\!$Y (crosses) and erg$_2\!$Y (squares), respectively,
for three sample paths.}
\label{figure1}       
\end{figure}

The results are presented in Table \ref{table1} and Figure \ref{figure1}. 

\begin{table}
\begin{center}
\caption{The marginal CDF (limiting loss function), its approximations and the corresponding relative errors (RE) for different values of $\epsilon$ and for three
realisations.
The number of time steps is $N=4\hspace{-1pt}\cdot\hspace{-1pt}10^{5}$ and the relative statistical error is $0.15$\%.}\label{table1}
\begin{tabularx}{\textwidth}{@{}YYYYYYYYYYYY@{}}
\addlinespace[-5pt]
  \toprule[.1em]
	\centering $\!\!\!\!\!\epsilon\quad$ & \textbf{\!\!\!\!\!\!Loss} & \textbf{appY} & \centering \textbf{RE}(\%) & \textbf{erg$_1$\!Y} & \centering \textbf{RE}(\%) & \textbf{erg$_2$\!Y} & \centering \textbf{RE}(\%) & \textbf{erg$_1$\!YZ} & \centering \textbf{RE}(\%) & \textbf{erg$_2$\!YZ} & \centering \textbf{RE}(\%) \tabularnewline
  \midrule
  	$\!\!\!\!\!\!10^{0}$ & $\!\!\!\!\!\!0.49715$ & $0.42303$ & $14.91$ & $0.41542$ & $16.44$ & $0.43350$ & $12.80$ & $0.51886$ & $4.37$ & $0.54713$ & $10.05$ \\[1pt]
	$\!\!\!\!\!\!10^{-1}$ & $\!\!\!\!\!\!0.53123$ & $0.54974$ & $3.48$ & $0.55471$ & $4.42$ & $0.57538$ & $8.31$ & $0.51886$ & $2.33$ & $0.54713$ & $2.99$ \\[1pt]
	$\!\!\!\!\!\!10^{-2}$ & $\!\!\!\!\!\!0.67626$ & $0.67967$ & $0.50$ & $0.69624$ & $2.95$ & $0.71826$ & $6.21$ & $0.51886$ & $23.27$ & $0.54713$ & $19.09$ \\[1pt]
	$\!\!\!\!\!\!10^{-3}$ & $\!\!\!\!\!\!0.58330$ & $0.58041$ & $0.50$ & $0.58835$ & $0.87$ & $0.61049$ & $4.66$ & $0.51886$ & $11.05$ & $0.54713$ & $6.20$ \\[1pt]
	$\!\!\!\!\!\!10^{-4}$ & $\!\!\!\!\!\!0.47947$ & $0.48005$ & $0.12$ & $0.47805$ & $0.30$ & $0.49815$ & $3.90$ & $0.51886$ & $8.22$ & $0.54713$ & $14.11$ \\
  \midrule
	$\!\!\!\!\!\!10^{0}$ & $\!\!\!\!\!\!0.45688$ & $0.39515$ & $13.51$ & $0.38492$ & $15.75$ & $0.40177$ & $12.06$ & $0.41543$ & $9.07$ & $0.43981$ & $3.74$ \\[1pt]
	$\!\!\!\!\!\!10^{-1}$ & $\!\!\!\!\!\!0.46628$ & $0.47136$ & $1.09$ & $0.46849$ & $0.47$ & $0.48729$ & $4.51$ & $0.41543$ & $10.91$ & $0.43981$ & $5.68$ \\[1pt]
	$\!\!\!\!\!\!10^{-2}$ & $\!\!\!\!\!\!0.48804$ & $0.48556$ & $0.51$ & $0.48412$ & $0.80$ & $0.50405$ & $3.28$ & $0.41543$ & $14.88$ & $0.43981$ & $9.88$ \\[1pt]
	$\!\!\!\!\!\!10^{-3}$ & $\!\!\!\!\!\!0.32785$ & $0.32987$ & $0.62$ & $0.31403$ & $4.22$ & $0.32807$ & $0.07$ & $0.41543$ & $26.71$ & $0.43981$ & $34.15$ \\[1pt]
	$\!\!\!\!\!\!10^{-4}$ & $\!\!\!\!\!\!0.41799$ & $0.42107$ & $0.74$ & $0.41327$ & $1.13$ & $0.43136$ & $3.20$ & $0.41543$ & $0.61$ & $0.43981$ & $5.22$ \\
	\midrule
	$\!\!\!\!\!\!10^{0}$ & $\!\!\!\!\!\!0.17821$ & $0.15721$ & $11.78$ & $0.13417$ & $24.71$ & $0.13932$ & $21.82$ & $0.15003$ & $15.81$ & $0.15806$ & $11.30$ \\[1pt]
	$\!\!\!\!\!\!10^{-1}$ & $\!\!\!\!\!\!0.13892$ & $0.15375$ & $10.68$ & $0.13075$ & $5.88$ & $0.13572$ & $2.30$ & $0.15003$ & $8.00$ & $0.15806$ & $13.78$ \\[1pt]
	$\!\!\!\!\!\!10^{-2}$ & $\!\!\!\!\!\!0.17593$ & $0.18456$ & $4.90$ & $0.16152$ & $8.19$ & $0.16803$ & $4.49$ & $0.15003$ & $14.72$ & $0.15806$ & $10.16$ \\[1pt]
	$\!\!\!\!\!\!10^{-3}$ & $\!\!\!\!\!\!0.25130$ & $0.25282$ & $0.61$ & $0.23195$ & $7.70$ & $0.24213$ & $3.65$ & $0.15003$ & $40.30$ & $0.15806$ & $37.10$ \\[1pt]
	$\!\!\!\!\!\!10^{-4}$ & $\!\!\!\!\!\!0.31246$ & $0.31196$ & $0.16$ & $0.29477$ & $5.66$ & $0.30801$ & $1.42$ & $0.15003$ & $51.99$ & $0.15806$ & $49.41$ \\[1pt]
	\bottomrule[.1em]
  \addlinespace[3pt]
\end{tabularx}
\end{center}
\end{table}

The approximate conditional CDF (appY), derived from Lyapunov's central limit theorem,
provides a good fit to the true conditional CDF for small $\epsilon$ across all samples.
Figure \ref{figure1} allows a comparison of this error to an asymptotic behaviour of order $\epsilon^{1/2}$, however, due to the irregular behaviour of the individual path realisations, no definitive conclusions are possible.

For linear averaging of the volatility function in $Y^{1}$, erg$_1\!$Y, which gives the correct conditional (on $W^x, W^y$) expectation of $X^1$ for $\epsilon\rightarrow 0$,
the error broadly decreases for decreasing $\epsilon$ but is significantly larger than in the CLT-based approximation.
A similar behaviour is observed for quadratic averaging erg$_2\!$Y, which gives the correct asymptotic second unconditional moment of $X^1$.

The approximations based on full (Y and Z) linear and quadratic averages, erg$_1\!$YZ and  erg$_2\!$YZ, respectively, are independent of $\epsilon$. It is seen from the last columns of Table \ref{table1} that they give a poor approximation to the true loss and are therefore not included in Figure \ref{figure1} for clarity. 

The latter observations are in line with \cite[Theorem 2.4]{Hambly:2020} who derive a limiting particle system for $\epsilon\to 0$ in which the  averaged squared volatility and a modified correlation coefficient appear. In the SPDE for the limit empirical measure, this is replaced by a linear average and a yet different correlation coefficient. As per \cite[Corollary 2.8]{Hambly:2020}, this indicates that  except for $\rho_y=0$, convergence is generally only observed in a distributional sense and not strongly.



\section{Weak approximation of loss function}\label{sec:weak}

In this section, we give an application to basket credit derivatives and analyse numerically
the accuracy of the approximations.


Let $a\in[0,1]$ and consider a call option on the limiting loss function at time $T$,
\begin{equation}\label{eq5.1}
C_{a} = \EE\left[(L_{T}-a)^{+}\right],
\end{equation}
with fixed default level $B$. For convenience, we assume that $\rho_{x}>0$. This type of payoffs is common in credit derivatives, e.g., in single tranche CDOs \cite{Bush:2011,Giles:2012}.

We can compute the call price via \eqref{eq4.18} by estimating the limiting loss function at time $T$ by Monte Carlo sampling of $W^{y,1}$, and then the outer expectation in \eqref{eq5.1} by sampling of $W^{x}$ and $W^{y}$ (abbreviated \emph{limCall}). Alternatively, for a large number of firms, we can approximate the call price by
\[
C_{N_{\hspace{-1.5pt}f},\hspace{1pt}a} = \EE\left[(L_{N_{\hspace{-1.5pt}f},T}-a)^{+}\right]
\]
and then estimate the expectation by a Monte Carlo average over discrete trajectories of $W^{x}$, $W^{y}$, $W^{x,1}, \hdots, W^{x,N_{\hspace{-1.5pt}f}}$, $W^{y,1}, \hdots, W^{y,N_{\hspace{-1.5pt}f}}$ (abbreviated \emph{firmsCall}). The latter method does not require that we simulate an inner expectation, which can be very expensive, but we lose the smoothness in the loss function. 

\subsection{Call price approximation by conditional CLT argument}\label{subsec:appCall}

Here, we approximate the call price using the approximate marginal CDF from \eqref{eq4.7.9}.
We can decompose the Brownian motion $W^{x}$ as
$W^{x} = \rho_{xy}W^{y} + \sqrt{1-\rho_{xy}^{2}}\hspace{1pt}\tilde{W}^{x}$,
where $W^{y}$ and $\tilde{W}^{x}$ are independent Brownian motions. Let $\mathcal{F}^{y}$ be the filtration generated by the market Brownian driver $W^{y}$. Conditional on the $\sigma$-algebra $\mathcal{F}^{y}_{T}$, 
\[
\int_{0}^{T}{\langle\sigma(\cdot+Z_{t})\rangle_{Y}\hspace{1pt}d\tilde{W}^{x}_{t}} \,\eqlaw\, \sqrt{\int_{0}^{T}{\langle\sigma(\cdot+Z_{t})\rangle^{2}_{Y}\hspace{1pt}dt}}\,W_{1},
\]
where $W_{1}$ is a standard normal random variable. Hence, the approximate conditional (on $\mathcal{F}^{y}_{T}$) law of the limiting loss function at time $T$ is that of
$\Phi\big(c_{Y,0}-c_{Y,1}W_{1}\big)$,
where
\begin{align}
c_{Y,0} &= \frac{B-\int_{0}^{T}{\langle\mu(\cdot+Z_{t})\rangle_{Y}\hspace{1pt}dt}-\rho_{x}\rho_{xy}\int_{0}^{T}{\langle\sigma(\cdot+Z_{t})\rangle_{Y}\hspace{1pt}dW^{y}_{t}}}{\sqrt{\int_{0}^{T}{\langle\sigma^{2}(\cdot+Z_{t})\rangle_{Y}\hspace{1pt}dt}-\rho_{x}^{2}\int_{0}^{T}{\langle\sigma(\cdot+Z_{t})\rangle_{Y}^{2}\hspace{1pt}dt}}}, \\
\label{eq5.8}
c_{Y,1} &= \rho_{x}\sqrt{1-\rho_{xy}^{2}}\hspace{1pt}\sqrt{\frac{\int_{0}^{T}{\langle\sigma(\cdot+Z_{t})\rangle^{2}_{Y}\hspace{1pt}dt}}{\int_{0}^{T}{\langle\sigma^{2}(\cdot+Z_{t})\rangle_{Y}\hspace{1pt}dt}-\rho_{x}^{2}\int_{0}^{T}{\langle\sigma(\cdot+Z_{t})\rangle^{2}_{Y}\hspace{1pt}dt}}}\hspace{1pt}.
\end{align}
Using a conditioning technique, we can express the call price as
\begin{equation}\label{eq5.9}
C_{a} = \EE\left[\EE\Big[(L_{T}-a)^{+} \,|\, \mathcal{F}^{y}_{T}\Big]\right].
\end{equation}
Upon noticing that $c_{Y,1}>0$, we can compute the inner expectation
\begin{align}\label{eq5.10}
\EE\Big[(L_{T}-a)^{+} \,|\, \mathcal{F}^{y}_{T}\Big] &\approx \EE\left[\Big(\Phi\big(c_{Y,0}-c_{Y,1}W_{1}\big)-a\Big)^{+} \,|\, \mathcal{F}^{y}_{T}\right] \nonumber\\[2pt]
&= \int_{-\infty}^{w_{0}}{\Big(\Phi\big(c_{Y,0}-c_{Y,1}w\big)-\Phi\big(\Phi^{-1}(a)\big)\Big)\phi(w)\hspace{1pt}dw} \nonumber\\[2pt]
&= \int_{-\infty}^{w_{0}}{\Phi\big(c_{Y,0}-c_{Y,1}w\big)\phi(w)\hspace{1pt}dw} - a\Phi(w_{0}),
\end{align}
with $\phi$ the standard normal PDF and
$w_{0} \!=\! \frac{c_{Y,0}-\Phi^{-1}(a)}{c_{Y,1}}\hspace{1pt}$.
By \cite{Owen:1980}, formula (10,010.1),
\begin{equation}\label{eq5.12}
\int_{-\infty}^{w_{0}}{\Phi\big(c_{Y,0}-c_{Y,1}w\big)\phi(w)\hspace{1pt}dw} = \bvn\left(\frac{c_{Y,0}}{\sqrt{1+c_{Y,1}^{2}}}\hspace{1pt},w_{0};\hspace{1pt}\frac{c_{Y,1}}{\sqrt{1+c_{Y,1}^{2}}}\right),
\end{equation}
where the bivariate normal CDF is
\begin{equation}\label{eq5.13}
\bvn\left(h,k;\rho\right) = \frac{1}{2\pi\sqrt{1-\rho^{2}}}\int_{-\infty}^{k}{\int_{-\infty}^{h}{\exp\left(-\frac{x^{2}-2\rho xy+y^{2}}{2(1-\rho^{2})}\right)\hspace{1pt}dx}\hspace{1pt}dy}.
\end{equation}
Combining \eqref{eq5.10}--\eqref{eq5.12} yields
\begin{eqnarray}\label{eq5.14}
\EE\Big[(L_{T}-a)^{+} \,|\, \mathcal{F}^{y}_{T}\Big] &\approx& \\
&&\hspace{-3 cm} \bvn\left(\frac{c_{Y,0}}{\sqrt{1+c_{Y,1}^{2}}}\hspace{1pt},\hspace{1pt}\frac{c_{Y,0}-\Phi^{-1}(a)}{c_{Y,1}}\hspace{1pt};\hspace{1pt}\frac{c_{Y,1}}{\sqrt{1+c_{Y,1}^{2}}}\right) - a\Phi\left(\frac{c_{Y,0}-\Phi^{-1}(a)}{c_{Y,1}}\right).
\end{eqnarray}
Finally, we discretize the two coefficients, i.e., $c_{Y,0}\approx \bar{c}_{Y,0}$ and $c_{Y,1}\approx \bar{c}_{Y,1}$, where
\begin{eqnarray}\label{eq5.15}
\bar{c}_{Y,0} &=& \frac{Bm^{-1}e^{-\frac{\xi^{2}}{k}(1-\rho_{y}^{2})}+\frac{1}{2}\hspace{1pt}me^{\frac{\xi^{2}}{k}(1-\rho_{y}^{2})}
\mathcal{I} -\rho_{x}\rho_{xy}e^{-\frac{\xi^{2}}{2k}(1-\rho_{y}^{2})}  \mathcal{M} }{\sqrt{\left(1-\rho_{x}^{2}e^{-\frac{\xi^{2}}{k}(1-\rho_{y}^{2})}\right)\delta t\sum_{n=0}^{N-1}{e^{2z_{t_{n}}}}}}, \ \textrm{with} \\
&& \hspace{-1 cm} \mathcal{I} = \sum_{n=0}^{N-1}{e^{2z_{t_{n}}}} \delta t, \;\;
\mathcal{M} = \sum_{n=0}^{N-1}{e^{z_{t_{n}}}\hspace{1pt}\Big(W^{y}_{t_{n+1}}-W^{y}_{t_{n}}\Big)},
\;\; and \;\;
\bar{c}_{Y,1} = \frac{\rho_{x}\sqrt{1-\rho_{xy}^{2}}}{\sqrt{e^{\frac{\xi^{2}}{k}(1-\rho_{y}^{2})}-\rho_{x}^{2}}}\hspace{1pt},
\nonumber
\end{eqnarray}
and estimate the outer expectation in \eqref{eq5.9} by Monte Carlo sampling of $W^{y}$.

\subsection{Call price approximation by ergodic averages}\label{subsec:ergCall}

\textit{$Y$-averages.} 
First, we approximate the call price by employing a linear or quadratic ergodic $Y^{1}$ average. Recall from \eqref{eq4.10} and \eqref{eq4.11} that
\begin{equation}\label{eq5.16.1}
L_{T} \approx \Phi\left(\frac{B-\int_{0}^{T}{\langle\mu(\cdot+Z_{t})\rangle_{Y}\hspace{1pt}dt}-\rho_{x}\int_{0}^{T}{\langle\sigma^{2-\lambda}(\cdot+Z_{t})\rangle^{\frac{1}{{2-\lambda}}}_{Y}\hspace{1pt}dW^{x}_{t}}}{\sqrt{(1-\rho_{x}^{2})\int_{0}^{T}{\langle\sigma^{2}(\cdot+Z_{t})\rangle_{Y}\hspace{1pt}dt}}}\right).
\end{equation}
Proceeding as before, we deduce that
\begin{eqnarray}\label{eq5.16.2}
\EE\Big[(L_{T}-a)^{+} \,|\, \mathcal{F}^{y}_{T}\Big] &\approx& \bvn\left(\frac{c_{Y,2}}{\sqrt{1+c_{Y,3}^{2}}}\hspace{1pt},\hspace{1pt}\frac{c_{Y,2}-\Phi^{-1}(a)}{c_{Y,3}}\hspace{1pt};\hspace{1pt}\frac{c_{Y,3}}{\sqrt{1+c_{Y,3}^{2}}}\right) \\
&& - a\Phi\left(\frac{c_{Y,2}-\Phi^{-1}(a)}{c_{Y,3}}\right), \nonumber
\end{eqnarray}
where
\begin{eqnarray}\label{eq5.16.3}
c_{Y,2} &=& \frac{B-\int_{0}^{T}{\langle\mu(\cdot+Z_{t})\rangle_{Y}\hspace{1pt}dt}-\rho_{x}\rho_{xy}\int_{0}^{T}{\langle\sigma^{2-\lambda}(\cdot+Z_{t})\rangle^{\frac{1}{{2-\lambda}}}_{Y}\hspace{1pt}dW^{y}_{t}}}{\sqrt{(1-\rho_{x}^{2})\int_{0}^{T}{\langle\sigma^{2}(\cdot+Z_{t})\rangle_{Y}\hspace{1pt}dt}}}, \\
\label{eq5.16.4}
c_{Y,3} &=& \frac{\rho_{x}\sqrt{1-\rho_{xy}^{2}}}{\sqrt{1-\rho_{x}^{2}}}\hspace{1pt}\sqrt{\frac{\int_{0}^{T}{\langle\sigma^{2-\lambda}(\cdot+Z_{t})\rangle^{\frac{2}{{2-\lambda}}}_{Y}\hspace{1pt}dt}}{\int_{0}^{T}{\langle\sigma^{2}(\cdot+Z_{t})\rangle_{Y}\hspace{1pt}dt}}}\hspace{1pt}.
\end{eqnarray}
As before, we discretize the two coefficients, i.e., $c_{Y,2}\approx \bar{c}_{Y,2}$ and $c_{Y,3}\approx \bar{c}_{Y,3}$, where
\begin{eqnarray}\label{eq5.16.5}
\bar{c}_{Y,2} &=& \frac{Bm^{-1}e^{-\frac{\xi^{2}}{k}(1-\rho_{y}^{2})}+\frac{1}{2}\hspace{1pt}me^{\frac{\xi^{2}}{k}(1-\rho_{y}^{2})} \mathcal{I}-
\rho_{x}\rho_{xy}e^{-\lambda\frac{\xi^{2}}{2k}(1-\rho_{y}^{2})} \mathcal{M}}{\sqrt{(1-\rho_{x}^{2})\delta t\sum_{n=0}^{N-1}{e^{2z_{t_{n}}}}}}, \\
\bar{c}_{Y,3} &=& \frac{\rho_{x}\sqrt{1-\rho_{xy}^{2}}}{\sqrt{1-\rho_{x}^{2}}}\hspace{1pt}e^{-\lambda\frac{\xi^{2}}{2k}(1-\rho_{y}^{2})},
\end{eqnarray}
and estimate the outer expectation in \eqref{eq5.9} by a 
sample average over $W^{y}$.

\textit{$Y$ and $Z$-averages.}  Last, we approximate the call price by linear and quadratic ergodic $Y^{1}$ and $Z$ average. 
Using \eqref{eq4.14} and \eqref{eq4.15}, we can deduce in a similar fashion
\begin{eqnarray}\label{eq5.18}
C_{a} &\approx& \bvn\left(\frac{c_{0}}{\sqrt{1+c_{1}^{2}}}\hspace{1pt},\hspace{1pt}\frac{c_{0}-\Phi^{-1}(a)}{c_{1}}\hspace{1pt};\hspace{1pt}\frac{c_{1}}{\sqrt{1+c_{1}^{2}}}\right) - a\Phi\left(\frac{c_{0}-\Phi^{-1}(a)}{c_{1}}\right), \\
\label{eq5.19}
\textrm{where} && c_{0} = \frac{Bm^{-1}e^{-\frac{\xi^{2}}{k}}+\frac{1}{2}\hspace{1pt}me^{\frac{\xi^{2}}{k}}T}{\sqrt{(1-\rho_{x}^{2})T}}, \qquad
c_{1} = \frac{\rho_{x}}{\sqrt{1-\rho_{x}^{2}}}\hspace{1pt}e^{-\lambda\frac{\xi^{2}}{2k}}.
\end{eqnarray}

\subsection{Expected loss}\label{subsec:expLoss}

In the special case of a linear payoff ($a=0$), the call price is simply the expected limiting loss function at time $T$ (abbreviated \emph{expLoss}). Using a conditioning technique and \eqref{eq4.3}, we can write the expected loss as
\begin{equation}\label{eq5.21}
\EE\left[L_{T}\right] = \EE\left[\EE\left[\Ind_{X^{1}_{T}\leq B} \,|\, \mathcal{F}^{x,y}_{T}\right]\right] = \Prob\Big(X^{1}_{T}\leq B\Big) = \EE\left[\Prob\Big(X^{1}_{T}\leq B \,|\, \mathcal{F}^{y,y_{1}}_{T}\Big)\right].
\end{equation}
From \eqref{eq1.1} 
we deduce that
\begin{equation}\label{eq5.22}
\EE\left[L_{T}\right] = \EE\left[\Phi\left(\frac{B-\int_{0}^{T}{\mu(Y^{1}_{t}+Z_{t})\hspace{1pt}dt}-\rho_{x}\rho_{xy}\int_{0}^{T}{\sigma(Y^{1}_{t}+Z_{t})\hspace{1pt}dW^{y}_{t}}}{\sqrt{(1-\rho_{x}^{2}\rho_{xy}^{2})\int_{0}^{T}{\sigma^{2}(Y^{1}_{t}+Z_{t})\hspace{1pt}dt}}}\right)\right],
\end{equation}
which can be estimated by a Monte Carlo average over samples of $W^{y}$ and $W^{y,1}$. Hence, this provides a much faster method 
in the special case of a linear payoff.

\subsection{Numerical tests}\label{subsec:numericsWeak}

We perform numerical tests for the weak errors with the different approximations. We fix the time horizon $T=1$ and the default level $B=-0.1$ as in Section \ref{sec:pathwise}, and assign the same values to the underlying model parameters as in \eqref{eq4.21}. Furthermore, we fix $\epsilon=4\hspace{-1pt}\cdot\hspace{-1.5pt}10^{-3}$ as in \cite{Dobson:2015}, a choice which corresponds to a mean-reversion time of 1.5 days, as observed from S\&P500 data (see \cite{Fouque:2003b}).

\begin{table}
\begin{center}
\caption{The call option price, its approximations and the corresponding relative errors (RE) for 3 different strikes $a\in\{0.00,0.05,0.10\}$. The number of time steps is $N=10^{4}$, the number of firms is $N_{\hspace{-1.5pt}f}\in\{5,150,100\}$ -- each value corresponds to one of the 3 strikes -- and the relative statistical errors are $0.15$\% for the call price and $0.05$\% for the approximations and the expected loss.}\label{table2}
\begin{tabularx}{\textwidth}{@{}LRRRRRR@{}}
	\addlinespace[-5pt]
  \toprule[.1em]
	& \multicolumn{2}{c}{\textbf{Strike} $\mathbf{=0.00}$} & \multicolumn{2}{c}{\textbf{Strike} $\mathbf{=0.05}$} & \multicolumn{2}{c}{\textbf{Strike} $\mathbf{=0.10}$} \\
	\textbf{Method} & \centering \hspace{2em}\textbf{Price} & \centering \hspace{2em}\textbf{RE} & \centering \hspace{2em}\textbf{Price} & \centering \hspace{2em}\textbf{RE} & \centering \hspace{2em}\textbf{Price} & \centering \hspace{2em}\textbf{RE} \tabularnewline
  \midrule
	\textbf{expLoss} & $0.18835$ & -- & -- & -- & -- & -- \\[1pt]
	\textbf{firmsCall} & $0.18843$ & $0.05$\% & $0.16170$ & -- & $0.14132$ & -- \\
	\midrule
	\textbf{appY} & $0.18878$ & $0.23$\% & $0.16155$ & $0.09$\% & $0.14078$ & $0.38$\% \\[1pt]
	\textbf{erg$_1$\!Y} & $0.18390$ & $2.36$\% & $0.15860$ & $1.92$\% & $0.13941$ & $1.35$\% \\[1pt]
	\textbf{erg$_2$\!Y} & $0.18872$ & $0.20$\% & $0.16342$ & $1.06$\% & $0.14410$ & $1.97$\% \\[1pt]
	\textbf{erg$_1$\!YZ} & $0.18262$ & $3.04$\% & $0.15724$ & $2.76$\% & $0.13789$ & $2.42$\% \\[1pt]
	\textbf{erg$_2$\!YZ} & $0.18912$ & $0.41$\% & $0.16376$ & $1.27$\% & $0.14423$ & $2.06$\% \\[1pt]
	\bottomrule[.1em]
	\addlinespace[3pt]
\end{tabularx}
\end{center}
\end{table}

The number of samples for the outer expectations in, e.g., \eqref{eq5.9},
was $1.2 \cdot 10^6$ and gave a small statistical error, estimated as the corrected sample standard deviation of the estimator divided by the value itself,
of $0.15$\% for the call price and $0.05$\% for the approximations and the expected loss.
We have verified numerically that the errors associated with the number of time steps $N$ and the number of firms $N_{\hspace{-1.5pt}f}$ from Table \ref{table2} match the statistical errors.

We infer from the data in Table \ref{table2} that $\epsilon=4\hspace{-1pt}\cdot\hspace{-1.5pt}10^{-3}$ gives a very small appY-approximation error throughout. 
Squared averaging (conditional on $Z$ or unconditional), where the first two moments of the $X$ process are matched, results in a very good approximation for a linear payoff, but in a worse approximation than the linear average for a non-linear payoff. 

\section{Conclusions}

It has recently been shown theoretically in \cite{Hambly:2020} that large pool models of processes with fast mean-reverting stochastic volatility
may be approximated by one-dimensional models with constant, averaged model parameters.
The limit as the mean-reversion speed goes to infinity is generally only attained in a distributional, but not in a strong sense.

We show in this paper how such averaged equations can be implemented numerically, but also observe that the approximation quality
is poor in both the strong and the weak sense in cases of interest.

The main finding of the paper is an improved approximation obtained by a central limit theorem argument, which leads to consistently good
accuracy both in a path-wise sense conditional on common noise, and in a weak sense when considering expected nonlinear functionals of the solution.
A theoretically rigorous analysis of this empirically improved estimator will be the topic of future research.


\end{document}